\documentclass[draft]{amsart}

\title[Truncations determining $p$-divisible groups up to
isogeny]{Traverso's isogeny conjecture for $p$-divisible groups}
\author{Marc-Hubert Nicole and Adrian Vasiu}

\date{\today; to appear in Rend. Sem. Mat. Univ. Padova}

\usepackage{amssymb}
\usepackage{graphicx, epsfig}
\usepackage{
latexsym, amsfonts, amscd, amsmath}

\newtheorem{thm}{Theorem}[section]
\newtheorem{prop}[thm]{Proposition}
\newtheorem{lem}[thm]{Lemma}
\newtheorem{cor}[thm]{Corollary}

\newtheorem{con}[thm]{Conjecture}

\newtheorem{rmk}[thm]{Remark}
\newtheorem{dfn}[thm]{Definition}
\newtheorem{exa}[thm]{Example}

\newtheorem{fact}[thm]{Fact}

\input amssym.def

\newcommand{\calD}{{\mathcal{D}}}

\newcommand{\calH}{{\mathcal{H}}}

\newcommand{\calN}{{\mathcal{N}}}

\newcommand{\End}{{\operatorname{End}}}

\newcommand{\Ker}{{\operatorname{Ker}}}

\newcommand{\Spec}{{\operatorname{Spec }}}
\newcommand{\Spf}{{\operatorname{Spf }}}

\def\db#1{{\fam\msbfam\relax#1}}

 \def\dbF{{\db F}}

 \def\dbN{{\db N}}
 
\def\dbQ{{\db Q}}

 \def\dbZ{{\db Z}}

\def\Q{\mathbb{Q}}

\def\Ker{\text{Ker}}

\def\End{\text{End}}
\def\Spec{\text{Spec}}

\def\leaderfill{\leaders\hbox to 1em
     {\hss.\hss}\hfill}

\def\finishproclaim{\par\rm
     \ifdim\lastskip\medskipamount\removelastskip
     \penalty55\medskip\fi}
\def\endproof{$\hfill \square$}
\def\proof{\par\noindent {\it Proof:}\enspace}

\def\Ref[#1]{\par\hang\indent\llap{\hbox to\parindent
     {[#1]\hfil\enspace}}\ignorespaces}
\def\Item#1{\par\smallskip\hang\indent\llap{\hbox to\parindent
     {#1\hfill$\,\,$}}\ignorespaces}
\def\ItemItem#1{\par\indent\hangindent2\parindent
     \hbox to \parindent{#1\hfill\enspace}\ignorespaces}

\def\Le{{\mathchoice{\,{\scriptstyle\le}\,}
  {\,{\scriptstyle\le}\,}
  {\,{\scriptscriptstyle\le}\,}{\,{\scriptscriptstyle\le}\,}}}

\def\arrowsim{\,\smash{\mathop{\to}\limits^{\lower1.5pt
  \hbox{$\scriptstyle\sim$}}}\,}

\def\doublemaprights#1#2#3#4{\raise3pt\hbox{$\mathop{\,\,\hbox to
#1pt{\rightarrowfill}\kern-30pt\lower3.95pt\hbox to
     #2pt{\rightarrowfill}\,\,}\limits_{#3}^{#4}$}}

\def\rightcapdownarrow{\raise9pt\hbox{$\ssize\cap$}\kern-7.75pt
     \Big\downarrow}

\def\rcapmapdown#1{\rightcapdownarrow\kern-1.0pt\vcenter{
     \hbox{$\scriptstyle#1$}}}

\def\rmapdown#1{\Big\downarrow\kern-1.0pt\vcenter{
     \hbox{$\scriptstyle#1$}}}
\def\rightsubsetarrow#1{{\ssize\subset}\kern-4.5pt\lower2.85pt
     \hbox to #1pt{\rightarrowfill}}
\def\longtwoheadedrightarrow#1{\raise2.2pt\hbox to #1pt{\hrulefill}
     \!\!\!\twoheadrightarrow}

\usepackage[usenames]{color}
\definecolor{Indigo}{rgb}{0.2,0.1,0.7}
\definecolor{Violet}{rgb}{0.5,0.1,0.7}
\definecolor{White}{rgb}{1,1,1}
\definecolor{Green}{rgb}{0.1,0.9,0.2}

\begin{document}

\maketitle \centerline{\it To Carlo Traverso, for his 60th
birthday}

\bigskip\noindent
{\bf ABSTRACT.} Let $k$ be an algebraically closed field of
characteristic $p>0$. Let $c,d\in\dbN$. Let $b_{c,d}\ge 1$ be the smallest
integer such that for any two $p$-divisible groups $H$ and $H^\prime$ over
$k$ of codimension $c$ and dimension $d$ the following assertion holds: If
$H[p^{b_{c,d}}]$ and $H^\prime[p^{b_{c,d}}]$ are isomorphic, then $H$ and
$H^\prime$ are isogenous. We show that $b_{c,d}=\lceil{cd\over
{c+d}}\rceil$.
This proves Traverso's isogeny conjecture for $p$-divisible groups over
$k$.


\section{Introduction}

Let $p\in\dbN$ be a prime. Let $k$ be an algebraically closed
field of characteristic $p$. Let $c,d\in\dbN$. Let
$H$ be a $p$-divisible group over $k$ of codimension $c$ and
dimension $d$; its height is $c+d$. It is well-known
(see [Ma], [Tr1, Thm. 3], [Tr2, Thm. 1], [Va, Cor. 1.3], and [Oo3, Cor.
1.7]) that there exists a minimal number $n_H\in\dbN$ such that $H$ is
uniquely
determined up to isomorphism by $H[p^{n_H}]$ (i.e., if $H^\prime$
is a $p$-divisible group over $k$ such that $H^\prime[p^{n_H}]$ is
isomorphic to $H[p^{n_H}]$, then $H^\prime$ is isomorphic to $H$).
For instance, we have $n_H\le cd+1$ (cf. [Tr1, Thm. 3]). This implies that
there exists a minimal number $b_H\in\dbN$ such that
the {\it isogeny class} of $H$ is uniquely determined by
$H[p^{b_H}]$. We call $b_H$ the {\it isogeny cutoff} of $H$.
We have $1\le b_H\le n_H$. Traverso speculated the following (cf. [Tr3,
\S40, Conj. 5]):

\begin{con} \label{traversocon} The isogeny cutoff $b_H$ is bounded
from above by $\lceil{cd\over {c+d}}\rceil$ i.e., we have $b_H\Le
\lceil{cd\over {c+d}}\rceil$.
\end{con}

\indent
The goal of this paper is to prove an optimal variant of the Conjecture:

\begin{thm} \label{traversoisogeny} Let $c,d\in\dbN$ and let $b_{c,d}\ge
1$ be the smallest integer such that for any two $p$-divisible groups $H$
and $H^\prime$ over $k$ of codimension $c$ and dimension $d$ the following
assertion holds: If $H[p^{b_{c,d}}]$ and $H^\prime[p^{b_{c,d}}]$ are
isomorphic, then $H$ and $H^\prime$ are isogenous. Then
$b_{c,d}=\lceil{cd\over {c+d}}\rceil$ and therefore we have $b_H\le
\lceil{cd\over
{c+d}}\rceil$.
\end{thm}

\noindent \indent In Section 2, we introduce notations and basic
invariants which pertain to $p$-divisible groups over $k$ and
which allow us to obtain a practical upper bound of $b_H$ (see
Corollary 2.13). In Section 3, we prove Theorem
\ref{traversoisogeny}; see Example 2.10 for a simpler proof in the
particular case when $H$ is isosimple. If the $a$-number of $H$ is
at most $1$, then the inequality $b_H\le \lceil{cd\over
{c+d}}\rceil$ is in essence due to Traverso (cf. [Tr3, Thm. of
\S21, and \S40]; see Theorem 3.1). We refer to Theorem 2.2 and
Example 3.2 for concrete examples with $b_H=1$ and with
$b_H=\lceil{cd\over {c+d}}\rceil$ (respectively).

\section{Estimates of the isogeny cutoff $b_H$}

Let $r:=c+d$ and $j:=\lceil{cd\over {r}}\rceil\in\dbN$. Let $W(k)$ be the
ring of Witt vectors with coefficients in $k$.
Let $\sigma$ be the Frobenius automorphism of $W(k)$ induced from
$k$. Let $(M,\phi)$ be the (contravariant) Dieudonn\'e module of
$H$. Thus $M$ is a free $W(k)$-module of rank $r$ and $\phi:M\to
M$ is a $\sigma$-linear endomorphism such that we have
$pM\subseteq \phi(M)$. We have $c=\dim_k(\phi(M)/pM)$ and
$d=\dim_k(M/\phi(M))$. Let $\vartheta:M\to M$ be the
$\sigma^{-1}$-linear Verschiebung map of
$\phi$; we have $\phi\vartheta=\vartheta\phi=p1_M$.

The Dieudonn\'e--Manin classification of $F$-isocrystals
over $k$ (see [Di, Thms. 1 and 2], [Ma, Ch. II, \S4], [De, Ch. IV]) states
that:

\medskip
$\bullet$ there exists $m\in\dbN$ such that we have a direct sum
decomposition $\break (M[{1\over p}],\phi)=\oplus_{i=1}^m
(M_i,\phi)$ into simple $F$-isocrystals over $k$, and

\smallskip
$\bullet$ there exist numbers $c_i,d_i\in\dbN\cup\{0\}$ which satisfy the
inequality $r_i:=d_i+c_i>0$, which are relative prime (i.e.,
$g.c.d.(c_i,d_i)=1$), and for which there exists a $B(k)$-basis for $M_i$
formed by elements fixed by $p^{-d_i}\phi^{r_i}$.

\medskip
\indent The unique slope of $(M_i,\phi)$ is $\alpha_i:={{d_i}\over
{r_i}}\in [0,1]\cap\dbQ$. To fix ideas, we assume that
$\alpha_1\le\alpha_2\le\cdots\le\alpha_m$. Each
$i\in\{1,\ldots,m\}$ gives a slope $\alpha_i$ with multiplicity
$r_i$. The Newton polygon $\calN_H$ of $H$ is a continuous,
piecewise linear, upward convex function $\calN_H:[0,r]\to [0,d]$
which for all $i\in\{1,\ldots,m\}$ has slope $\alpha_i$ on the
interval $[\sum_{\ell=1}^{i-1} r_{\ell},\sum_{\ell=1}^i r_{\ell}]$; we have $\calN_H(0)=0$ and $\calN_H(r)=d$.
As the field $k$ is algebraically closed, the isogeny
class of $H$ is uniquely determined by the Newton polygon
$\calN_H$. The finite set $\calN_{c,d}$ of Newton polygons we thus
obtain by varying $H$, can be partially ordered as follows: we say
that $\calN_1$ lies above (resp. strictly above) $\calN_2$ if for
all $t\in [0,r]$ we have $\calN_1(t)\ge\calN_2(t)$ (resp. if for
all $t\in [0,r]$ we have $\calN_1(t)\ge\calN_2(t)$ and moreover
$\calN_1\neq\calN_2$). This partial order is convenient when
studying the variation of the Newton polygon in families. We will
use the notation $\calN_{\bullet}$ to denote the Newton polygon of
$\bullet$, where $\bullet$ is a $p$-divisible group over a field
that contains $k$.

\noindent Let $H_{c_i,d_i}$ be the $p$-divisible group over $k$
whose Dieudonn\'e module $(M_{c_i,d_i},\phi_{c_i,d_i})$ has the
property that there exists a $W(k)$-basis
$\{e_0,\ldots,e_{r_i-1}\}$ for $M_{c_i,d_i}$ such that we have
$\phi_{c_i,d_i}(e_l)=e_{d_i+l}$ if $l\in\{0,\ldots,c_i-1\}$ and
$\phi_{c_i,d_i}(e_l)=pe_{d_i+l}$ if $l\in\{c_i,\ldots,r_i-1\}$
(here the lower right indices of $e$ are taken modulo $r_i$). It
is well-known that $H_{c_i,d_i}$ has the following two properties
(see [dJO, Lem. 5.4]): (i) it has codimension $c_i$, dimension
$d_i$, and unique Newton polygon slope $\alpha_i$, and (ii)  its
endomorphism ring $\End(H_{c_i,d_i})$ is the maximal order in the
simple $\dbQ_p$-algebra $\End(H_{c_i,d_i})\otimes_{\dbZ_p} \Q_p$
of invariant ${\alpha_i}\in \dbQ/\dbZ$. In addition, these two
properties determine $H_{c_i,d_i}$ up to isomorphism (see [Ma] and
\cite[\S 1.1]{Oo4}). \noindent Let $H_0:= \prod_{i=1}^{m}
H_{c_i,d_i}$; we have $\calN_{H_0}=\calN_H$. The Dieudonn\'e
module of $H_0$ is $(M_0,\phi_0):=\prod_{i=1}^{m}
(M_{c_i,d_i},\phi_{c_i,d_i})$.

\begin{dfn} \emph{(Oort)} \label{wellknown}
We say that a $p$-divisible group $H$ over $k$ is {\it minimal} if
it is isomorphic to the $p$-divisible group $H_0$ determined
uniquely from the Newton polygon $\calN_H$ of $H$.
\end{dfn}

\noindent The $p$-divisible group $H_0$ is the \emph{unique} (up
to isomorphism) minimal $p$-divisible over $k$ that is isogenous
to $H$.

\begin{thm} \label{thm1}
A minimal $p$-divisible group $H$ over $k$ is determined up to
isomorphism by $H[p]$ (thus we have $b_H=n_H=1$).
\end{thm}

\proof For the isoclinic case (i.e., when the pairs $(c_i,d_i)$ do
not depend on $i$) we refer to [Va1, Exa. 3.3.6]. For the general
case we refer to either  {\cite[Thm. 1.2]{Oo4}} or [Va2, Thm.
1.8].\endproof

\begin{dfn} \label{minimalheight}
By the {\it minimal height} $q_H$ of a $p$-divisible group $H$
over $k$ we mean the smallest non-negative integer $q_H$ such that
there exists an isogeny $\psi_0: H\to H_0$ to a minimal
$p$-divisible group, whose kernel $\Ker(\psi_0)$ is annihilated by
$p^{q_H}$.
\end{dfn}

\begin{dfn} \label{anumber}
By the {\it a-number} $a_H$ of a $p$-divisible group $H$ over $k$
we mean the number
$\dim_k(\pmb{\alpha}_p,H[p])=\dim_k(M/\phi(M)+\vartheta(M))\in\dbN\cup\{0\}$,
where $\pmb{\alpha}_p$ is the local-local group scheme of order
$p$.
\end{dfn}

See
[Oo2, Prop. 2.8] for a proof of the following specialization
trick.

\begin{prop} \label{onlygoodthinkofoort} For every $p$-divisible group $H$
over $k$, there exists a $p$-divisible group $\calH$ over $k[[x]]$ that
has the following two properties:

\medskip
{\bf (i)} its fibre over $k$ is $H$;

\smallskip
{\bf (ii)} if $\overline{k((x))}$ is an algebraic closure of
$k((x))$,  then $\calH_{\overline{k((x))}}$ has the same
Newton polygon as $H$ and its $a$-number is at most $1$.
\end{prop}

\begin{dfn} \label{weakisogenybound}
By the {\it weak isogeny cutoff} of a $p$-divisible group $H$ over
$k$ we mean the smallest number $\tilde b_H\in\dbN$ such that the
following two properties hold:

\medskip
{\bf (i)} if $H^\prime$ is a $p$-divisible group over $k$ such that
$H^\prime[p^{\tilde b_H}]$ is isomorphic to $H[p^{\tilde b_H}]$, then its
Newton polygon $\calN_{H^\prime}$ is not strictly above $\calN_H$;

\smallskip
{\bf (ii)} there exists a $p$-divisible group $\calH$ over
$k[[x]]$ that has the following three properties:

\medskip
{\bf (ii.a)} its fibre over $k$ is $H$;

\smallskip
{\bf (ii.b)} the fibre $\calH_{k((x))}$ over $k((x))$ has the
same Newton polygon as $H$;

\smallskip
{\bf (ii.c)} the isogeny cutoff of
$\calH_{\overline{k((x))}}$ is at most $\tilde b_H$ and the
$a$-number of $\calH_{\overline{k((x))}}$ is at most $1$.
\end{dfn}

\noindent Note that property (i) holds for any level greater or
equal to $b_H$. Thus the existence of the number $\tilde
b_H\in\dbN$ is implied by Proposition 2.5.

\begin{fact} \label{fact1} If $H^{\text{t}}$ is the Cartier dual of $H$,
then $q_H=q_{H^{\text{t}}}$, $b_H=b_{H^{\text{t}}}$, and $\tilde
b_H=\tilde b_{H^{\text{t}}}$.
\end{fact}

\proof
This follows from Cartier duality (the Cartier dual of $H_{c_i,d_i}$ is
$H_{d_i,c_i}$).
\endproof

\begin{lem} \label{lemma1} The isogeny cutoff $b_H$ is the smallest
natural number such that for every element $g\in \pmb{GL}_M(W(k))$
congruent to $1_M$ modulo $p^{b_H}$, the Dieudonn\'e module
$(M,g\phi)$ is isogenous to $(M,\phi)$.
\end{lem}

\proof Let $t\in\dbN$. The arguments proving either [Va1, Cor.
3.2.3] or \cite[Thm. 2.2 (a)]{NV} show that: (i) if $g\in
\pmb{GL}_M(W(k))$ is congruent to $1_M$ modulo $p^t$, then every
$p$-divisible group $H_g$ over $k$ whose Dieudonn\'e module is
isomorphic to $(M,g\phi)$, has the property that $H_g[p^t]$ is
isomorphic to $H[p^t]$, and (ii) if a $p$-divisible group $H^\prime$
over $k$ has the property that $H^\prime[p^t]$ is isomorphic to
$H[p^t]$, then there exists $g\in \pmb{GL}_M(W(k))$ congruent to
$1_M$ modulo $p^t$ and such that the Dieudonn\'e module of $H^\prime$
is isomorphic to $(M,g\phi)$. The Lemma follows from these two
statements and the very definition of the isogeny cutoff $b_H$.\endproof

\begin{lem} \label{good} For every $p$-divisible group $H$ over $k$, the
isogeny cutoff $b_H$ is bounded from above by the minimal height $q_H$
plus
$1$ i.e., we have $b_H\le q_H+1$.
\end{lem}

\proof Let $\psi: H\to H_0$ be an isogeny whose kernel is
annihilated by $p^{q_H}$. Let $(M_0,\phi_0)\hookrightarrow
(M,\phi)$ be the monomorphism of Dieudonn\'e modules associated to the
isogeny $\psi$; we
will identify $M_0$ with its image under this monomorphism. We
have $p^{q_H}M\subseteq M_0\subseteq M$. Let $g\in
\pmb{GL}_M(W(k))$ be congruent to $1_M$ modulo $p^{q_H+1}$. We
write $g=1_M+p^{q_H+1}e$, where $e\in\End(M)$. We have
$p^{q_H+1}e(M_0)\subseteq p^{q_H+1}e(M)\subseteq
p^{q_H+1}M\subseteq pM_0$. Thus we have $g\in\pmb{GL}_{M_0}(W(k))$
and moreover $g$ is congruent  to $1_{M_0}$ modulo $p$. Therefore
the reductions modulo $p$ of the two triples
$(M_0,\phi,\vartheta)$ and $(M_0,g\phi,\vartheta g^{-1})$
coincide. As the Dieudonn\'e module $(M_0,\phi)$ is minimal (being
isomorphic to $(M_0,\phi_0)$), from Theorem 2.2 we get that it is
isomorphic to $(M_0,g\phi)$. The isogeny class of $(M,g\phi)$
(resp. of $(M,\phi)$) is the same as of $(M_0,g\phi)$ (resp. of
$(M_0,\phi)$). From the last two sentences we get that the
Dieudonn\'e modules $(M,g\phi)$ and $(M,\phi)$ are isogenous. From
this and Lemma 2.8, we get that $b_H\le q_H+1$.
\endproof

\begin{exa} \label{dJO}
\emph{Suppose that $m=1$ i.e., $H$ is an isosimple $p$-divisible
group. Let $\theta_0:H_0\to H_0$ be an endomorphism such that
$\End(H_0)=W(\dbF_{p^r})[\theta_0]$ and $\theta_0^r=p$, cf. [Ma,
Ch. III, \S4, 1] and [dJO, Lem. 5.4]. From [dJO, 5.8] we get that
there exist inclusions $(M_0,\phi_0)\hookrightarrow
(M,\phi)\hookrightarrow (\theta_0^{-(c-1)(d-1)}M_0,\phi_0)$
between Dieudonn\'e modules over $k$. Let $\psi_0:H\to H_0$ be the
isogeny defined by the first inclusion
$(M_0,\phi_0)\hookrightarrow (M,\phi)$. Its kernel is annihilated
by $\theta_0^{(c-1)(d-1)}=p^{{{(c-1)(d-1)}\over r}}$ and
therefore we have $q_H\le \lceil{{(c-1)(d-1)}\over{r}}\rceil=j-1$.
Thus $b_H\le j$, cf. Lemma 2.9.}
\end{exa}

\begin{rmk} \label{induction}
\emph{The function $f:(0,\infty)\times (0,\infty)\to (0,\infty)$
defined by the rule $f(x,y)={{xy}\over {x+y}}$ is subadditive
i.e., for all $x_1,x_2,y_1,y_2\in (0,\infty)$, we have an
inequality $f(x_1,y_1)+f(x_2,y_2)\le f(x_1+x_2,y_1+y_2)$. But the
function $g(x,y):=\lceil{f(x,y)}\rceil$ is not subadditive. Due to
this and the plus $1$ part of Lemma 2.9, our attempts to use Lemma
2.9 in order to prove Theorem 1.2 by induction on $m\in\dbN$,
failed. On the other hand, in most cases Lemma 2.9 provides better
upper bounds of $b_H$ than those provided by the following
Proposition.}
\end{rmk}

\begin{prop} \label{better} For every $p$-divisible group $H$ over $k$,
the
isogeny cutoff $b_H$ is bounded from above by the weak isogeny cutoff
$\tilde b_H$ i.e., we have $b_H\le\tilde b_H$.
\end{prop}

\proof Let $\overline{k((x))}$ be an algebraic closure of
$k((x))$. Let $\calH$ be a $p$-divisible group over $k[[x]]$ of
constant Newton polygon such that its fibre over $k$ is
$H$ and the isogeny cutoff $b$ of
$\calH_{\overline{k((x))}}$ is at most $\tilde b_H$, cf. property
(ii) of Definition 2.6. We recall that the statement that $\calH$
has constant Newton polygon means that the Newton
polygons of $H$ and $\calH_{\overline{k((x))}}$ coincide i.e.,
$\calN_H =\calN_{\calH_{\overline{k((x))}}}$. For the proof of this
Proposition it is irrelevant what the $a$-number of
$\calH_{\overline{k((x))}}$ is.

Let $H^\prime$ be a $p$-divisible group over $k$ such that
$H^\prime[p^{\tilde b_H}]=H[p^{\tilde b_H}]$. Based on [Il, Thm.
4.4. f)] we get that for all $n\in\dbN$, there exists a
$p$-divisible group $\calH^\prime_n$ over $\Spec(k[[x]]/(x^n))$
that lifts $H^\prime$ and such that $\calH^\prime_n[p^{\tilde
b_H}]=\calH[p^{\tilde b_H}]\times_{k[[x]]} k[[x]]/(x^n)$. Based on
loc. cit., we can assume that $\calH^\prime_{n+1}$ restricted to
$\Spec(k[[x]]/(x^n))$ is $\calH^\prime_n$. Therefore the
$\calH^\prime_n$'s glue together to define a $p$-divisible group
$\calH^{\prime\text{f}}$ over the formal scheme $\Spf(k[[x]])$. We
recall that the categories of $p$-divisible groups over
$\Spf(k[[x]])$ and $\Spec(k[[x]])$ are naturally equivalent, cf.
[dJ, Lem. 2.4.4]. Thus let $\calH^\prime$ be the $p$-divisible
group over $\Spec(k[[x]])$ that corresponds naturally to
$\calH^{\prime\text{f}}$.

\noindent The $p$-divisible group $\calH^\prime$ lifts $H^\prime$
and we have
$$\calH^\prime[p^{\tilde b_H}]=\text{proj.}\text{lim.}_{n\to\infty}
\calH^\prime_n[p^{\tilde b_H}]=\calH[p^{\tilde b_H}].$$
As $\calH^\prime[p^{\tilde b_H}]=\calH[p^{\tilde b_H}]$ and
$b\le\tilde b_H$, from the very definition of $b$, we get that
$\calH^\prime_{\overline{k((x))}}$ has the same Newton polygon as
$\calH_{\overline{k((x))}}$; thus
$\calN_{\calH^\prime_{\overline{k((x))}}} = \calN_H$. As the
Newton polygons go up under specialization (see [De, Ch. IV, \S7, Thm.]),
we conclude that the Newton
polygon $\calN_{H^\prime}$ of $H^\prime$ is above the Newton polygon
$\calN_H$ of $H$. But $\calN_{H^\prime}$ is not strictly above
$\calN_H$, cf. property (i) of Definition 2.6. From the last two
sentences we get that $\calN_{H^\prime}=\calN_H$. This implies that
$b_H\le\tilde b_H$.\endproof

\medskip
From Lemma 2.10 and Proposition 2.12 we get:

\begin{cor} \label{practicalupperbound} For every $p$-divisible group $H$
over $k$, we
have the following inequality $b_H\le\min\{\tilde b_H,q_H+1\}$.
\end{cor}

\section{The proof of Theorem 1.2}

In this Section we prove Theorem 1.2. We begin by proving the
following particular case of Conjecture 1.1 which in essence is due
to Traverso (cf. [Tr3, Thm. of \S21, and \S40]).

\begin{thm} \label{thmTraverso}
Suppose that $a_H\le 1$. Then we have $b_H\le j.$
\end{thm}

\proof We first recall how to compute the Newton polygon $\calN_H$
of $H$ (cf. [Ma], [De, Ch. IV, Lem. 2, pp. 82--84], [Tr3, \S21,
Thm.], and [Oo1, Lem. 2.6]). Let $v_p:W(k)\to\dbN\cup\{0,\infty\}$ be the normalized valuation of $W(k)$; thus $v_p(p)=1$.
Let $x\in M$ be such that its reduction modulo $p$ generates the
$k$-vector space $M/\phi(M)+\vartheta(M)$. Let
$a_0,a_1,\ldots,a_c,b_1,\ldots,b_d\in W(k)$ be such that the map
$$\psi:=\sum_{i=0}^c a_{c-i}\phi^i+\sum_{\ell=1}^d
b_{\ell}\vartheta^{\ell}:M\to M$$
annihilates $x$. It is easy to see that we can choose $x$ such that
$v_p(a_0)=v_p(b_d)=0$ and
$\{x,\phi(x),\ldots,\phi^{c-1}(x),\vartheta(x),\ldots,\vartheta^d(x)\}$ is
a $W(k)$-basis for $M$. For $\ell\in\{1,\ldots,d\}$ let
$a_{c+\ell}:=p^{\ell}b_{\ell}$; thus $a_i\in W(k)$ is well defined for all
$i\in\{0,\ldots,r\}$. Let
$$Q_x(t):=t^r+\sum_{i \in \{1, \ldots, r\}, a_i \neq 0}
v_p(a_i)t^{r-i}=t^r+\sum_{i \in \{1, \ldots, r\}, a_i \neq 0}
v_p(\sigma^d(a_i))t^{r-i}\in\dbZ[t].$$ We recall that the Newton
polygon of $Q_x(t)$ (more precisely, of the $(r+1)$-tuple
$(a_0,\ldots,a_r)$) is the greatest continuous, piecewise linear,
upward convex function $\calN_x:[0,r]\to [0,d]$ with the property
that for all $i\in\{0,\ldots,r\}$ we have $\calN_x(i)\le
v_p(a_i)$. Then we have
$$\calN_H=\calN_x.\leqno (1)$$
\indent
To check Formula (1) we view $M$ as a $W(k)[F]$-module,
where $F \cdot \lambda = \lambda^{\sigma} F$ and where $F$ acts on $M$ as
$\phi$ does. The $W(k)[F]$-module $M$ is isogenous to
the $W(k)[F]$-module $M^\prime:=W(k)[F]/W(k)[F]\psi^\prime$, where
$$\psi^\prime= F^d \psi:= \sum_{i=0}^r \sigma^d(a_i)F^{r-i}.$$
But the Newton polygon of the $W(k)[F]$-module $M^\prime$ is $\calN_x$,
cf. [De, Ch. IV, Lem. 2, pp. 82--84]. Thus Formula (1) holds.

We recall that $j = \lceil{{cd}\over{r}}\rceil$. As $c,d>0$, we
have $j\ge 1$. Let $g\in\pmb{GL}_M(W(k))$ be congruent to $1_M$
modulo $p^j$. Let $H_g$ be a $p$-divisible group over $k$ whose
Dieudonn\'e module is isomorphic to $(M,g\phi)$. As $j\ge 1$,
$H_g[p]$ is isomorphic to $H[p]$ and therefore $a_{H_g}=a_H\le 1$.
Based on Lemma 2.8, to prove the Theorem it suffices to show that
the Newton polygons $\calN_{H_g}$ and $\calN_H$ coincide. By replacing
$(\phi,\vartheta)$ with $(g\phi,\vartheta g^{-1})$, the map
$\psi:M\to M$ which annihilates $x$ gets replaced by another  map
$$\psi_g:=\sum_{i=0}^c a_{g,c-i}\phi^i+\sum_{\ell=1}^d
b_{g,\ell}\vartheta^{\ell}:M\to M$$
which annihilates $x$, where $a_{g,c-i}$ is congruent to $a_{c-i}$
modulo $p^j$ and where $b_{g,\ell}$ is congruent to $b_{\ell}$
modulo $p^j$. For $\ell\in\{1,\ldots,d\}$ let
$a_{g,c+\ell}:=p^{\ell}b_{g,\ell}$; thus $a_{g,i}\in W(k)$ is well defined
for all $i\in\{0,\ldots,r\}$ and it is congruent to $a_i$ modulo
$p^{j+\min\{0,i-c\}}$. The pair $(Q_x(t),\calN_x)$ gets replaced
by the pair $(Q_{g,x}(t),\calN_{g,x})$, where
$Q_{g,x}(t):=t^r+\sum_{i\in\{1,\ldots,r\},a_i\neq 0} v_p(a_{g,i})t^{r-i}$ and where
$\calN_{g,x}$ is the Newton polygon of $Q_{g,x}(t)$. As $\calN_x$
is upward convex and as $\calN_x(0)=0$ and $\calN_x(r)=d$, we have
$\calN_x(t)\le {{dt}\over r}$ for all $t\in [0,r]$. Thus
$j=\lceil{{cd}\over {r}}\rceil\ge\lceil{\calN_x(c)}\rceil\ge
\calN_x(c)\ge\calN_x(i)$ for all
$i\in\{0,\ldots,c\}$. As $r>d$ and $j\ge {{cd}\over r}$, for
$i\in\{c+1,\ldots,r\}$ we have
$j+i-c\ge {{di}\over r}$. From the last two sentences we get that
$j+\min\{0,i-c\}\ge {{di}\over r}\ge\calN_x(i)$ for all
$i\in\{0,\ldots,r\}$. From these inequalities and the fact that for all
$i\in\{0,\ldots,r\}$ the elements $a_{g,i}$ and $a_i$ are
congruent modulo $p^{j+\min\{0,i-c\}}$, we easily get that
$\calN_{g,x}=\calN_x$. From this identity and Formula (1) we get that
$\calN_{H_g}=\calN_H$.\endproof

\subsection {End of the proof of Theorem 1.2.}
Based on Proposition 2.12, to prove the inequality $b_H\le j$
it suffices to show that $\tilde b_H\le j$. We recall that $\calN_{c,d}$
is the
set of Newton polygons of $p$-divisible groups over $k$ of
codimension $c$ and dimension $d$. Let $\calD_H$ be the class of
$p$-divisible groups of codimension $c$ and dimension $d$ over $k$
whose Newton polygons are strictly above $\calN_H$.

We prove the
inequality $\tilde b_H\le j$ by decreasing induction on
$\calN_H\in\calN_{c,d}$. Thus we can assume that for every
$\bullet\in\calD_H$ we have $\tilde b_{\bullet}\le j$; as
$b_{\bullet}\le \tilde b_{\bullet}$ (cf. Proposition 2.12) we also
have $b_{\bullet}\le j$. To show that $\tilde b_H\le j$, let $\calH$ be a
$p$-divisible group over $k[[x]]$ such that the properties (ii.a) to
(ii.c) of Definition
2.6 hold. Let $b$ be the isogeny cutoff of
$\calH_{\overline{k((x))}}$. From the very definition of $\tilde b_H$ we
get that
$$\tilde b_H\le\max\{b,b_{\bullet}|\bullet\in\calD_H\}.$$
As $b\le j$ (cf. Theorem 3.1 applied to $\calH_{\overline{k((x))}}$) and
as $b_{\bullet}\le j$
for all $\bullet\in\calD_H$ (cf. the inductive assumption), we have
$\tilde b_H\le j$. This ends the induction.

Thus $\tilde b_H\le j$ and
therefore $b_H\le j$. As $H$ is an arbitrary $p$-divisible group over $k$
of codimension $c$ and dimension $d$, we get that $b_{c,d}\le j$. If
$j=1$, then we obviously have $b_{c,d}=1$. If $j\ge 2$, then the below
Example shows that there exist $p$-divisible groups $H$ over $k$ of
codimension $c$ and dimension $d$ and such that we have $b_H=j$. This
implies that for $j\ge 2$, we have $b_{c,d}\ge j$. We conclude that for
all values of $j\in\dbN$ we have $b_{c,d}=j=\lceil{cd\over {r}}\rceil$.
\endproof

\begin{exa}
\emph{Suppose that
$j\ge 2$. Let $s\in\dbN$ be the smallest number such that $r$
divides $cd-s$; we have $j-1={{cd-s}\over r}$. We assume that
there exists an element $x\in M$ such that the $r$-tuple
$(e_1,\ldots,e_r):=(x,\phi(x),\ldots,\phi^{c-1}(x),\vartheta^d(x),\ldots,\vartheta(x))$
is an ordered $W(k)$-basis for $M$ and we have an equality
$\vartheta^d(x)=\phi^c(x)$. Thus $\phi(e_i)=e_{i+1}$ if
$i\in\{1,\ldots,c\}$ and  $\phi(e_i)=pe_{i+1}$ if
$i\in\{c+1,\ldots,r\}$. We have $\phi^r(e_i)=p^de_i$ and therefore
all Newton polygon slopes of $H$ are ${d\over r}$; thus
$m=g.c.d.(c,r)$. Replacing the equality $\vartheta^d(x)=\phi^c(x)$
by the equality $\tilde\vartheta^d(x)=\tilde\phi^c(x)-p^{j-1}x$,
we get a $p$-divisible group $\tilde H$ over $k$ whose Dieudonn\'e
module $(M,\tilde\phi)$ is such that $\tilde\phi(e_i)=e_{i+1}$ if
$i\in\{1,\ldots,c-1\}$, $\tilde\phi(e_c)=p^{j-1}e_1+e_{c+1}$, and
$\tilde\phi(e_i)=pe_{i+1}$ if $i\in\{c+1,\ldots,r\}$.  As
$(\tilde\phi,\tilde\vartheta)$ and $(\phi,\vartheta)$ are
congruent modulo $p^{j-1}$, $\tilde H[p^{j-1}]$ is isomorphic to
$H[p^{j-1}]$. We have $\calN_H(t)={{dt}\over r}$ but the piecewise
linear function $\calN_{\tilde H}(t)$ changes slope at $c$; more
precisely we have $\calN_{\tilde H}(c)=j-1$ (cf. Formula (1)
applied to $\tilde H$, to $x\in M$, to the function
$\tilde\psi:=-p^{j-1}\tilde\phi^0+\tilde\phi^c-\tilde\vartheta^d:M\to M$,
and to the polynomial $\tilde
Q_x(t):=t^r+(j-1)t^d+d$). Thus the Newton polygon slopes of
$\tilde H$ are ${{j-1}\over c}$ (with multiplicity $c$) and
$1-{{j-1}\over d}$ (with multiplicity $d$) and therefore they are
different from ${d\over r}$. This implies that $b_H\ge j$ and
therefore (cf. the inequality $b_H\le j$ proved in the previous paragraph)
we have $b_H=j$.}
\end{exa}

\begin{exa} \label{supersingular}
\emph{Suppose that $c=d$ and that $H$ is as in Example 3.2. Then
$q_H={{c-1}\over 2}$ (cf. [NV, Rmk. 3.3]) and $j=\lceil{c\over
2}\rceil$. If $c$ is odd, then $b_H=j={{c+1}\over 2}=q_H+1$; moreover, for
$\tilde H:=(\dbQ_p/\dbZ_p)\oplus H\oplus\mu_{p^{\infty}}$ we also have
$b_{\tilde H}=j={{c+1}\over 2}=q_{\tilde H}+1$. Thus for $c=d$,
the inequality $b_H\le q_H+1$ (see Lemma 2.9) is optimal in general. }
\end{exa}

\medskip\noindent
{\bf Acknowledgments.}

\medskip
The first author has been supported by the Japanese Society for
the Promotion of Science (JSPS PDF) while working on this paper at
the University of Tokyo. The second author would like to thank University
of Arizona for good conditions with which to write this note.

\bigskip
\hbox{Marc-Hubert Nicole,\;\;\;Email: nicole@ms.u-tokyo.ac.jp}
\hbox{Address: University of Tokyo, Department of Mathematical
Sciences,} \hbox{Komaba, 153-8914, Tokyo, Japan.}

\medskip

\hbox{Adrian Vasiu,\;\;\;Email: adrian@math.arizona.edu}
\hbox{Address: University of Arizona, Department of Mathematics,
617 N. Santa Rita Ave.,} \hbox{P.O. Box 210089, Tucson, AZ 85721-0089,
U.S.A.}

\end{document}